\documentclass[12pt]{article}
\usepackage{amsmath}
\usepackage{amssymb}

\def\digamma{{\cal F}}

\newtheorem{theorem}{Theorem}

\newtheorem{proposition}[theorem]{Proposition}
\newtheorem{remark}[theorem]{Remark}

\newenvironment{proof}[1][Proof]{\textbf{#1.} }{\ \rule{0.5em}{0.5em}}
\input{tcilatex}

\begin{document}

\title{\textbf{Finite-Time Blowup and Existence of Global Positive Solutions
of a Semi-Linear SPDE }}
\author{\textsc{Marco Dozzi} \and \textsc{Jos\'{e} Alfredo L\'{o}pez-Mimbela}%
}
\date{ }
\maketitle

\begin{abstract}
\noindent We consider stochastic equations of the prototype $du(t,x) =\left(
\Delta u(t,x)+u(t,x)^{1+\beta}\right)dt+\kappa u(t,x)\,dW_{t}$ on a smooth
domain $D\subset \mathord{\rm I\mkern-3.6mu R\:\!\!}^d$, with Dirichlet
boundary condition, where $\beta$, $\kappa$ are positive constants and $\{W_t
$, $t\ge0\}$ is a one-dimensional standard Wiener process. We estimate the
probability of finite time blowup of positive solutions, as well as the
probability of existence of non-trivial positive global solutions.

 \bigskip

\end{abstract}

\section{Introduction}

\label{section1}

Let $D\subset \mathord{\rm I\mkern-3.6mu R\:\!\!}^{d}$ be a bounded domain
with smooth boundary $\partial D.$ We consider a semilinear equation of the
form 
\begin{eqnarray}  \label{E1}
du(t,x) &=&\left( \Delta u(t,x)+G(u(t,x))\right) dt+\kappa
u(t,x)\,dW_{t},\quad t>0,  \notag \\
u(0,x) &=&f(x)\geqq 0,\quad x\in D, \\
u(t,x) &=&0,\quad t\geqq 0,\quad x\in \partial D,  \notag
\end{eqnarray}
where $G:\mathord{\rm I\mkern-3.6mu R\:\!\!}\rightarrow 
\mathord{\rm
I\mkern-3.6mu R\:\!\!}_{+}$ is locally Lipschitz and satisfies 
\begin{equation}  \label{G}
G(z)\geq Cz^{1+\beta }\quad \mbox{for all } \ z>0,
\end{equation}
$C,$ $\beta $ and $\kappa $ are given positive numbers, $\{W_{t},\ t\geq 0\}$
is a standard one-dimensional Brownian motion on a stochastic basis ($\Omega
,\mathcal{F} ,(\mathcal{F} _{t},t\geqq 0),P$), and $f:D\rightarrow 
\mathord{\rm I\mkern-3.6mu
R\:\!\!}_{+}$ is of class $C^{2}$ and not identically zero. We assume (\ref%
{G}) in sections 1 to 3 only; it is replaced by (\ref{G2}) in section 4.
Since we do not assume $G$ to be Lipschitz, blowup of the solution (\ref{E1}%
) in finite time can not be excluded, and our aim is to give estimates of
the probability of blowup and conditions for the existence of a global
solution of (\ref{E1}). A (random) time $T$ is called blowup time of $u$ if 
\begin{equation*}
\lim_{t\nearrow T}\sup_{x\in D}| u(t,x)| =+\infty \text{ \ \ }P-\mbox{a.s.} 
\text{ on }\{T<+\infty \}.
\end{equation*}
In the classical (deterministic) case where $G(z)=z^{1+\beta}$ and $\kappa=0$%
, it is well-known that for a nonnegative $f\in L^2(D)$, the condition  
\begin{equation}  \label{Dirichlet}
\int_Df(x)\psi(x)\,dx > \lambda_1^{1/\beta}
\end{equation}
already implies finite time blowup of \eqref{E1}. Here $\lambda_1>0$ is the
first eigenvalue of the Laplacian on $D$, and $\psi$ the corresponding
eigenfunction normalized so that $\|\psi\|_{L^1}=1$.

The existence, uniqueness and trajectorial regularity of global solutions of
parabolic equations perturbed by a time-homogenous white noise have been
investigated by different methods (see e.g. Chueshov and Vuillermot \cite%
{Chueshov and Vuillermot}, Denis et al. \cite{Denis et al}, Gy\"{o}ngy and
Rovira \cite{Gyongy and Rovira}, Krylov \cite{Krylov}, Lototski and
Rozovskii \cite{Lototski and Rozovskii}, Mikulevicius and Pragarauskas \cite%
{Mikulevicius and Pragarauskas}). Several types of solutions have been
proposed (see especially the last cited reference for strong solutions), and
the regularity results show that the solution is much smoother in the space
variable than for equations perturbed by space-dependent white noise.

Let us recall the notions of weak and mild solutions of (1) we are going to
use here. Let $\tau \leq +\infty $ be a stopping time. A continous $\mathcal{%
F}_{t}$-adapted random field $u=\{u(t,x),$ $t\geq 0,$ $x\in D\}$ is a \emph{%
weak solution} of (1) on the interval $]0,\tau \lbrack $ provided that, for
every $\varphi \in C^{2}(D)$ vanishing on $\partial D$, there holds 
\begin{eqnarray*}
\int_{D}u(t,x)\varphi (x)\,dx &=&\int_{D}f(x)\varphi
(x)\,dx+\int_{0}^{t}\int_{D}[u(s,x)\Delta \varphi (x)+G(u(s,x))\varphi
(x)]\,dx\,ds \\
&&+\kappa \int_{0}^{t}\int_{D}u(s,x)\varphi (x)\,dx\,dW_{s}\quad \mbox{$P$
-- a.s.}
\end{eqnarray*}%
for all $t\in \lbrack 0,\tau \lbrack $. Let $\{S_{t}$, $t\geqq 0\}$ be the
semigroup of $d$-dimensional Brownian motion killed at the boundary of $D$.
A continous $\mathcal{F}_{t}$-adapted random field $u=\{u(t,x),$ $t\geqq 0,$ 
$x\in D\}$ is a \emph{mild solution} of (1) on the interval $]0,\tau \lbrack 
$ if it satisfies 
\begin{equation*}
u(t,x)=S_{t}f(x)+\int_{0}^{t}\left[ S_{t-r}(G(u(r,\cdot ))(x)\,dr+\kappa
S_{t-r}(u(r,\cdot ))(x)\,dW_{r}\right] \quad \mbox{$P$-a.s. and $x$-a.e. in
$D$}
\end{equation*}%
for all $t\in ]0,\tau \lbrack $ (see e.g. \cite{Pazy}, Chapter IV). We refer
to \cite{Gyongy and Rovira} for background on existence of weak and mild
solutions, and for their equivalence under local Lipschitz conditions on $G$%
. Let us note that the results in \cite{Gyongy and Rovira} hold for a more
general class of second order differential operators which includes the
Laplacian as a special case. The positivity of the solution of (\ref{E1})
follows from comparison theorems (see e.g. Berg\'{e} et al. \cite{Berge} or
Manthey and Zausinger \cite{Manthey}).

Our aim in this communication is to study the blowup behaviour of $u$ by
means of a related random partial differential equation (see (\ref{Eq.3})
below). In section 3 we describe the blowup behaviour of the solution $v$ of
this random partial differential equation in terms of the first eigenvalue
and the first eigenfunction of the Laplace operator on $D$. This is done by
solving explicitly a stochastic equation in the time variable which is
obtained from the weak form of (\ref{Eq.3}). The solution of this
differential equation can be written in terms of integrals of exponential
Brownian motion with drift. The results of Dufresne \cite{Dufresne} and Yor 
\cite{Yor} on the law of these integrals easily imply estimates for the
probability of existence of a global solution, or of blowup in finite time
of $u$ and $v$. In section 4 sufficient conditions for $v$ to be a global
solution are given in terms of the semigroup of the Laplace operator using
recent sharp results on its transition density. These conditions show in
particular that the initial condition $f$ has to be small enough in order to
avoid for a given $G$ the blowup of $v$, and that the presence of noise may
help to prevent blowup. The results of section 4 can be used to investigate
the blowup behavior of $u$ by means of conditions (\ref{G}) and (\ref{G2}).

\section{A related random partial differential equation}

\label{section2}

In this section we investigate the random partial differential equation 
\begin{eqnarray}  \label{Eq.3}
\frac{\partial v}{\partial t}(t,x) &=&\Delta v(t,x)-\frac{\kappa^2}{2}%
v(t,x)+e^{-\kappa W_{t}}G(e^{\kappa W_{t}}v(t,x)),\quad t>0,\quad x\in D, 
\notag \\
v(0,x) &=&f(x),\quad x\in D, \\
v(t,x) &=&0,\text{ \ }x\in \partial D.  \notag
\end{eqnarray}
This equation is understood trajectorywise and classical results for partial
differential equations of parabolic type apply to show existence, uniqueness
and positivity of a solution up to eventual blowup (see e.g. Friedman \cite%
{Friedman} Chapter 7, Theorem 9).

\begin{proposition}
\label{proposition1}  Let $u$ be a weak solution of (\ref{E1}). Then the
function $v$ defined by 
\begin{equation*}
v(t,x)=e^{-\kappa W_{t}}u(t,x),\quad t\geq 0,\quad x\in D.
\end{equation*}
solves (\ref{Eq.3}).
\end{proposition}

\textbf{\noindent Remark } Proposition \ref{proposition1} implies in
particular that (\ref{E1}) possesses a strong local solution. 

\proof Recall that It\^{o}'s formula states that $\{e^{-\kappa W_{t}}$, $%
t\geq 0\}$ is the semimartingale given by 
\begin{equation*}
e^{-\kappa W_{t}}=1-\kappa \int_{0}^{t}e^{-\kappa W_{s}}dW_{s}+\frac{\kappa
^{2}}{2} \int_{0}^{t}e^{-\kappa W_{s}}\,ds.
\end{equation*}
Let us write $u(t,\varphi )\equiv \int_{D}u(t,x)\varphi (x)\,dx.$ Then a
weak solution of (\ref{E1}) can be written as 
\begin{equation*}
u(t,\varphi )=u(0,\varphi )+\int_{0}^{t}u(s,\Delta \varphi
)\,ds+\int_{0}^{t}G(u)(s,\varphi )\,ds+\kappa \int_{0}^{t}u(s,\varphi
)\,dW_{s}.
\end{equation*}
Therefore, for $\varphi $ fixed$,$ $\{u(t,\varphi )1_{[0,\tau [}(t),$ $t\geq
0\}$ is again a semimartingale. By applying the integration by parts formula
(see e.g. Klebaner \cite{Klebaner}, Ch. 8) we get 
\begin{eqnarray*}
v(t,\varphi ) &:=&\int_{D}v(t,x)\varphi (x)\,dx \\
&=&v(0,\varphi )+\int_{0}^{t}e^{-\kappa W_{s}}\,du(s,\varphi
)+\int_{0}^{t}u(s,\varphi )\left( -\kappa e^{-\kappa W_{s}}\,dW_{s}+\frac{%
\kappa ^{2}}{2}e^{-\kappa W_{s}}\,ds\right) \\
&&+\left[e^{-\kappa W_{\cdot}},u(\cdot ,\varphi )\right](t),
\end{eqnarray*}
where the quadratic variation is given by 
\begin{equation*}
\left[ e^{-\kappa W_{\cdot }},u(\cdot ,\varphi )\right](t)=-\int_{0}^{t}%
\kappa ^{2}e^{-\kappa W_{s}}u(s,\varphi )\,ds,\quad t\geq 0.
\end{equation*}
Therefore, 
\begin{eqnarray*}
v(t,\varphi ) &=&v(0,\varphi )+\int_{0}^{t}e^{-\kappa W_{s}}\left(
u(s,\Delta \varphi )+G(u)(s,\varphi )\right)\, ds+\kappa
\int_{0}^{t}e^{-\kappa W_{s}}u(s,\varphi )\,dW_{s} \\
&&-\kappa \int_{0}^{t}e^{-\kappa W_{s}}u(s,\varphi )\,dW_{s}+\frac{\kappa
^{2}}{2}\int_{0}^{t}e^{-\kappa W_{s}}u(s,\varphi )\,ds-\kappa
^{2}\int_{0}^{t}e^{-\kappa W_{s}}u(s,\varphi )\,ds \\
&=&v(0,\varphi )+\int_{0}^{t}\left[ v(s,\Delta \varphi )+e^{-\kappa
W_{s}}G(e^{\kappa W_{\cdot }}v)(s,\varphi )-\frac{\kappa ^{2}}{2}v(s,\varphi
)\right] \,ds.
\end{eqnarray*}
Moreover, by self-adjointness of the Laplacian, and the fact that $\varphi
(x)=0 $ for $x\in \partial D$, 
\begin{equation*}
v(s,\Delta \varphi )=\int_{D}v(s,x)\Delta \varphi (x),dx=\int_{D}\Delta
v(s,x)\varphi (x)\,dx=\Delta v(s,\varphi ).
\end{equation*}%
\hfill$\blacksquare$

\section{An estimate of the probability of blowup}

\label{section3}

Without loss of generality, let us assume that $C=1$ in (\ref{G}). Let $\psi 
$ be the eigenfunction corresponding to the first eigenvalue $\lambda_1 $ of
the Laplacian on $D$, normalized by $\int_{D}\psi (x)\,dx=1.$ It is
well-known that $\psi $ is strictly positive on $D$. Due to Proposition \ref%
{proposition1} we have that 
\begin{equation*}
v(t,\psi ) =v(0,\psi )+\int_{0}^{t}\left[v(s,\Delta \psi )\,-\frac{\kappa
^{2}}{2}v(s,\psi )\right]\,ds+\int_{0}^{t}e^{-\kappa W_{s}}G(e^{\kappa
W.}v)(s,\psi )\,ds. 
\end{equation*}
\bigskip Moreover, 
\begin{equation}
v(s,\Delta \psi )=-\lambda_1 v(s,\psi ),  \label{self-adjointness}
\end{equation}
and, due to (\ref{G}), 
\begin{equation}
\int_{D}e^{-\kappa W_{s}}G(e^{\kappa W_{s}}v(s,x))\psi (x)\,dx\ \geq\
e^{\kappa \beta W_{s}}\int_{D}v(s,x)^{1+\beta }\psi (x)\,dx .  \label{aux1}
\end{equation}
By Jensen's inequality 
\begin{equation}
\int_{D}v(s,x)^{1+\beta }\psi (x)\,dx\ \geq \ \left[ \int_{D}v(s,x)\psi
(x)\,dx \right] ^{1+\beta }\ =\ v(s,\psi )^{1+\beta },  \label{aux2}
\end{equation}
and therefore 
\begin{equation*}
\frac{d}{dt}v(t,\psi )\geq -\left( \lambda_1 +\frac{\kappa ^{2}}{2}\right)
v(t,\psi )+e^{\kappa \beta W_{t}}v(t,\psi )^{1+\beta }.
\end{equation*}
Hence $v(t,\psi )\ge I(t)$ for all $t\ge0$, where $I(\cdot)$ solves 
\begin{equation*}
\frac{d}{dt}I(t)=-\left( \lambda_1 +\frac{\kappa ^{2}}{2}\right)
I(t)\,+e^{\kappa \beta W_{s}}I(t)^{1+\beta }\,,\text{ \ }I(0)=v(0,\psi ),
\end{equation*}
and is given by 
\begin{equation*}
I(t)\ = \ e^{-(\lambda_1 +\kappa ^{2}/2)t}\left[ v(0,\psi )^{-\beta }-\beta
\int_{0}^{t}e^{-(\lambda_1 +\kappa ^{2}/2)\beta s+\kappa \beta W_{s}}\,ds %
\right] ^{-\frac{1}{\beta }},\quad 0\leq t<\tau ,
\end{equation*}
with 
\begin{equation}  \label{tautau}
\tau :=\inf \left\{ t\geq 0\left|\ \int_{0}^{t}e^{-(\lambda_1 +\kappa
^{2}/2)\beta s+\kappa \beta W_{s}}\,ds\geq \frac{1}{\beta }v(0,\psi
)^{-\beta }\right.\right\}.
\end{equation}
It follows that $I$ exhibits finite time blowup on the event $[\tau<\infty]$%
.  Since $I\leqq v(\cdot ,\psi ),$ $\tau $ is an upper bound for the blow-up
time of $v(\cdot ,\psi )$, and therefore for the blowup times of $v$ and $u$.

\begin{remark}
\label{Bandle-et-al} \emph{The same formula for the blow-up time of a {%
stochastic }differential equation, containing a stochastic integral with
respect to $W$, has been obtained in Bandle et al \cite{Bandle}. The
argument based on the first eigenvalue (and the corresponding eigenfunction)
of the Laplace operator on $D$ is applied there directly to $u$, and leads
to a stochastic differential inequality for $u(t,\psi )$. The associated
stochastic differential equation can again be solved explicitly by means of
the It\^{o} calculus, and the same formula as above is obtained for the
blowup time of the solution of this equation. Both approaches are therefore
equivalent, but the approach in \cite{Bandle} requires a more complicated
comparison theorem for stochastic differential inequalities. }
\end{remark}

Let us now give an estimate for the probability of blowup in finite time of $%
v$. From (\ref{tautau}), 
\begin{eqnarray}  \label{Dirichlet2}
P[\tau =+\infty ] &=&P\left[ \int_{0}^{t}\exp (-(\lambda_1 +\kappa
^{2}/2)\beta s+\kappa \beta W_{s})\,ds<\frac{1}{\beta }v(0,\psi )^{-\beta } 
\text{ for all }t>0\right]  \notag \\
&=&P\left[ \int_{0}^{\infty }\exp (-(\lambda_1 +\kappa ^{2}/2)\beta s+\kappa
\beta W_{s})\,ds\leqq \frac{1}{\beta }v(0,\psi )^{-\beta }\right] \\
&=&P\left[ \int_{0}^{\infty }\exp (2\widehat{\beta }W_{s}^{(\mu )})\,ds\leqq 
\frac{1}{\beta }v(0,\psi )^{-\beta }\right] ,  \notag
\end{eqnarray}
where $W_{s}^{(\mu )}:=\mu s+W_{s}$, $\mu :=-(\lambda_1 +\kappa
^{2}/2)/\kappa $, and $\hat{\beta}:=\kappa \beta /2$. Setting $\widehat{ \mu 
}=\mu /\widehat{\beta }$ we get 
\begin{equation}  \label{doubt}
P[\tau =+\infty ]=P\left[ \frac{4}{\kappa ^{2}\beta ^{2}}\int_{0}^{\infty
}\exp (2W_{s}^{(\widehat{\mu })})\,ds\leqq \frac{1}{\beta }v(0,\psi
)^{-\beta }\right] .
\end{equation}
It follows from \cite{Yor} (Chapter 6, Corollary 1.2) that 
\begin{equation*}
\int_{0}^{\infty }\exp (2W_{s}^{(\widehat{\mu })})\,ds\ =\ \frac{1}{ 2Z_{-%
\widehat{\mu }}}
\end{equation*}
in distribution, where $Z_{-\widehat{\mu }}$ is a random variable with law $%
\Gamma (-\widehat{ \mu }),$ i.e. $P(Z_{-\widehat{\mu }}\in dy)=\frac{1}{%
\Gamma (-\widehat{\mu }) }e^{-y}y^{-\widehat{\mu }-1}dy.$ We get therefore
(see also formula 1.10.4(1) in \cite{Borodin and Salminen (2002)}) 
\begin{equation*}
P[\tau =+\infty ]=\int_{0}^{\frac{1}{\beta }v(0,\psi )^{-\beta }}h(y)dy,
\end{equation*}
where 
\begin{equation*}
h(y)\ =\ \frac{(\kappa ^{2}\beta ^{2}y/2)^{(2\lambda_1 +\kappa ^{2})/\kappa
^{2}\beta }}{y\Gamma ((2\lambda_1 +\kappa ^{2})/(\kappa ^{2}\beta ))}\exp
\left( -\frac{2}{\kappa ^{2}\beta ^{2}y}\right). 
\end{equation*}
In this way we have proved the following

\begin{proposition}
\label{proposition2} The probability that the solution of (1) blows up in
finite time is lower bounded by $\int_{\frac{1}{\beta }v(0,\psi )^{-\beta
}}^{+\infty }h(y)\,dy.$
\end{proposition}

\begin{remark}
\label{Remarks}

\emph{\hskip-5pt\textbf{.1} Notice that formula 1.10.4(1) in \cite{Borodin
and Salminen (2002)} expresses the probability density function of $%
\int_{0}^{t}\exp (-(\lambda_1 +\kappa ^{2}/2)\beta s+\kappa \beta W_{s})\,ds$
in terms of the Kummer functions for\bigskip\ $\widehat{\mu }<2.$ }

\emph{\noindent \textbf{Remark 4.2} By putting $\kappa=0$ we get $v=u$ and,
moreover, in \eqref{Dirichlet2} we obtain that $P[\tau=+\infty]=0$ or $1$
according to $\int_Df(x)\psi(x)\,dx > \lambda_1^{1/\beta}$ or $%
\int_Df(x)\psi(x)\,dx \le \lambda_1^{1/\beta}$, which is a probabilistic
counterpart to \eqref{Dirichlet}. }
\end{remark}

\section{Non explosion of $v$}

We consider again equation (\ref{Eq.3}), but we assume now that $\kappa \neq
0$ and that $G:\mathord{\rm I\mkern-3.6mu R\:\!\!}_{+}\rightarrow 
\mathord{\rm I\mkern-3.6mu
R\:\!\!}_{+}$ satisfies $G(0)=0$, $G(z)/z$ is increasing and 
\begin{equation}
G(z)\leq \Lambda z^{1+\beta }\quad \mbox{for all }\ z>0,  \label{G2}
\end{equation}
where $\Lambda$ and $\beta$ are certain positive numbers. Let $\{S_t$, $%
t\geq 0\}$ again denote the semigroup of $d$-dimensional Brownian motion
killed at the boundary of $D$. Recall that Equation (\ref{Eq.3}) can be
re-written as 
\begin{equation}
v(t,x)=e^{-\kappa ^{2}t/2}S_{t}f(x)+\int_{0}^{t}e^{-\kappa
^{2}(t-r)/2}S_{t-r}\left( e^{-\kappa W_{r}}G\left(e^{\kappa W_{r}}v(r,\cdot
)\right)\right) (x)\,dr.  \label{mild}
\end{equation}
We give now a sufficient condition for the existence of a global solution of
(\ref{Eq.3}).

\begin{theorem}
\label{theorem3} Assume that $f$ satisfies 
\begin{equation}
\Lambda \beta \int_{0}^{\infty }e^{\kappa \beta W_{r}}\Vert e^{-\kappa
^{2}r/2}S_{r}f\Vert _{\infty }^{\beta }\,dr<1.  \label{cond1}
\end{equation}
Then Equation (\ref{Eq.3}) admits a global solution $v(t,x)$ that satisfies 
\begin{equation}
0\leq v(t,x)\leq \frac{e^{-\kappa ^{2}t/2}S_{t}f(x)}{\left( 1-\Lambda \beta
\int_{0}^{t}e^{\kappa \beta W_{r}}\Vert e^{-\kappa ^{2}r/2}S_{r}f\Vert
_{\infty }^{\beta }\,dr\right) ^{\frac{1}{\beta }}},\,\qquad t\geq 0.
\label{ResTh3}
\end{equation}
\end{theorem}

\proof Defining 
\begin{equation*}
B(t)=\left( 1-\Lambda \beta \int_{0}^{t}e^{\kappa \beta W_{r}}\Vert
e^{-\kappa ^{2}r/2}S_{r}f\Vert _{\infty }^{\beta }\,dr\right) ^{-\frac{1}{
\beta }},\quad t\geq 0,
\end{equation*}
we get $B(0)=1$ and 
\begin{equation*}
\frac{dB}{dt}(t)=\Lambda e^{\kappa \beta W_{t}}\Vert e^{-\kappa
^{2}t/2}S_{t}f\Vert _{\infty }^{\beta }B^{1+\beta}(t),
\end{equation*}
which implies 
\begin{equation*}
B(t)=1+\Lambda \int_{0}^{t}e^{\kappa \beta W_{r}}\Vert e^{-\kappa
^{2}r/2}S_{r}f\Vert _{\infty }^{\beta }B^{1+\beta }(r)\,dr.
\end{equation*}
Suppose now that $(t,x)\mapsto V_{t}(x)$ is a nonnegative continuous
function such that $V_{t}(\cdot )\in C_{0}(D)$, $t\geq 0$, and 
\begin{equation}
e^{-\kappa ^{2}t/2}S_{t}f(x)\leq V_{t}(x)\leq B(t)e^{-\kappa
^{2}t/2}S_{t}f(x),\quad t\geq 0,\quad x\in D.  \label{V}
\end{equation}
Let 
\begin{equation*}
R(V)(t,x):=e^{-\kappa ^{2}t/2}S_{t}f(x)+\int_{0}^{t}e^{-\kappa
W_{r}}e^{-\kappa ^{2}(t-r)/2}S_{t-r}\left( G(e^{\kappa W_{r}}V_{r}(\cdot
))\right) (x)\,dr.
\end{equation*}
Then, 
\begin{eqnarray}
\lefteqn{R(V)(t,x)}  \notag \\
&=&e^{-\kappa ^{2}t/2}S_{t}f(x)+\int_{0}^{t}e^{-\kappa W_{r}}e^{-\kappa
^{2}(t-r)/2}S_{t-r}\left( \frac{G(e^{\kappa W_{r}}V_{r}(\cdot ))}{
V_{r}(\cdot )}V_{r}(\cdot )\right) (x)\,dr  \notag \\
&\leq &e^{-\kappa ^{2}t/2}S_{t}f(x)+\int_{0}^{t}e^{-\kappa W_{r}}e^{-\kappa
^{2}(t-r)/2}S_{t-r}\left( \frac{G(e^{\kappa W_{r}}B(r)\Vert e^{-\kappa
^{2}r/2}S_{r}f\Vert _{\infty })}{B(r)\Vert e^{-\kappa ^{2}r/2}S_{r}f\Vert
_{\infty }}V(r)\right) (x)\,dr  \notag \\
&\leq &e^{-\kappa ^{2}t/2}S_{t}f(x)+\Lambda \int_{0}^{t}e^{\kappa \beta
W_{r}}B^{1+\beta }(r)\Vert e^{-\kappa ^{2}r/2}S_{r}f\Vert _{\infty }^{\beta
}e^{-\kappa ^{2}(t-r)/2}S_{t-r}(e^{-\kappa ^{2}r/2}S_{r}f)(x)\,dr  \notag \\
&=&e^{-\kappa ^{2}t/2}S_{t}f(x)\left[ 1+\Lambda \int_{0}^{t}e^{\kappa \beta
W_{r}}B^{1+\beta }(r)\Vert e^{-\kappa ^{2}r/2}S_{r}f\Vert _{\infty }^{\beta
}\,dr\right] =e^{-\kappa ^{2}t/2}S_{t}f(x)B(t),  \label{RV}
\end{eqnarray}
where to obtain the first inequality we used the rightmost inequality in (%
\ref{V}) and the fact that $G(z)/z$ is increasing, and to obtain the second
inequality we used (\ref{G2}). Consequently, 
\begin{equation*}
e^{-\kappa ^{2}t/2}S_{t}f(x)\leq R(V)(t,x)\leq B(t)e^{-\kappa
^{2}t/2}S_{t}f(x),\quad t\geq 0,\quad x\in D.
\end{equation*}
Let 
\begin{equation*}
v_{t}^{0}(x):=e^{-\kappa ^{2}t/2}S_{t}f(x)\quad \text{and \ \ }
v_{t}^{n+1}(x)=R(v^{n})(t,x),\quad n=0,1,2,\ldots .
\end{equation*}
Letting $n\rightarrow \infty $ yields, for $t\geq 0$ and $x\in D$, 
\begin{equation*}
0\leq v(t,x)= \lim_{n\rightarrow \infty }v_{t}^{n}(x)\leq B(t)e^{-\kappa
^{2}t/2}S_{t}f(x)\leq \frac{e^{-\kappa ^{2}t/2}S_{t}f(x)}{ \left( 1-\Lambda
\beta \int_{0}^{t}e^{\kappa \beta W_{r}}\Vert e^{-\kappa ^{2}r/2}S_{r}f\Vert
_{\infty }^{\beta }dr\right) ^{1/\beta }.}
\end{equation*}
Hence, $v(t,x)$ is a global solution of (\ref{mild}) due to the monotone
convergence theorem. \hfill $\blacksquare $

\bigskip

\noindent\textbf{Remark.} If we modify (\ref{Eq.3}) and (\ref{mild}) by
replacing $G(e^{\kappa W_{r}}v(t,x))$ by $G(v(t,x))$, then a global positive
solution still exists for all $f$ small enough, even if the inequality in (%
\ref{G2}) holds only for $z \in (0,C^{\ast })$, where $C^{\ast }$ is some
positive constant. In fact, if $f$ satisfies 
\begin{equation}
\Vert f\Vert _{\infty }\leq C^{\ast }\left( 1-\Lambda \beta \int_{0}^{\infty
}e^{-\kappa W_{r}}\Vert e^{-\kappa ^{2}r/2}S_{r}f\Vert _{\infty }^{\beta
}\,dr\right) ^{\frac{1}{\beta }},  \label{cond2}
\end{equation}
then Theorem \ref{theorem3} still holds if we replace the factor $e^{\kappa
\beta W_r}$ in (\ref{cond1}) and (\ref{ResTh3}) by the factor $e^{-\kappa
W_r}$. We only have to verify that assuming $z \in (0,C^{\ast })$ in (\ref%
{G2}) already implies the second inequality in (\ref{RV}): 
\begin{eqnarray*}
\Vert e^{-\kappa ^{2}t/2}S_{t}f\Vert _{\infty } &\leq &\Vert f\Vert _{\infty
} \\
&\leq &C^{\ast }\left( 1-\Lambda \beta \int_{0}^{t}e^{-\kappa W_{r}}\Vert
e^{-\kappa ^{2}r/2}S_{r}f\Vert _{\infty }^{\beta }\,dr\right) ^{ \frac{1}{%
\beta }}=\frac{C^{\ast }}{B^{\ast }(t)}\quad \text{for all }t\geqq 0,
\end{eqnarray*}
where 
\begin{equation*}
B^{\ast }(t)=\left( 1-\Lambda \beta \int_{0}^{t}e^{-\kappa W_{r}}\Vert
e^{-\kappa ^{2}r/2}S_{r}f\Vert _{\infty }^{\beta }\,dr\right) ^{-\frac{1}{
\beta }}.
\end{equation*}
This yields 
\begin{equation}
B^{\ast }(t)\Vert e^{-\kappa ^{2}t/2}S_{t}f\Vert _{\infty }\in (0,C^{\ast })
\end{equation}
for all $t\geq 0$, since $f\not\equiv 0$.

\bigskip

Let us now proceed to derive a sufficient condition for (\ref{cond1}) in
terms of the transition kernels $\{p_{t}(x,y)$, $t>0\}$ of $\{S_{t}$, $t\geq
0\}$ and the first eigenvalue $\lambda_1$ and corresponding eigenfunction $%
\psi$. We recall the following sharp bounds for $\{p_{t}(x,y)$, $t>0\}$,
which we borrowed from Ouhabaz and Wang \cite{Ouhabaz and Wang}.

\begin{theorem}
\label{Ouhabaz} Let $\psi >0$ be the first Dirichlet eigenfunction on a
connected bounded $C^{1,\alpha }$- domain in $%
\mathord{\rm
I\mkern-3.6mu R\:\!\!}^{d}$, where $\alpha >0$ and $d\geq 1$, and let $%
p_{t}(x,y)$ be the corresponding Dirichlet heat kernel. There exists a
constant $c>0$ such that, for any $t>0$, 
\begin{equation*}
\max \left\{ 1,\frac{1}{c}t^{-(d+2)/2}\right\} \leq e^{\lambda
_{1}t}\sup_{x,y}\frac{p_{t}(x,y)}{\psi (x)\psi (y)}\leq 1+c(1\wedge
t)^{-(d+2)/2}e^{-(\lambda _{2}-\lambda _{1})t},
\end{equation*}
where $\lambda _{2}>\lambda _{1}$ are the first two Dirichlet eigenvalues.
This estimate is sharp for both short and long times.
\end{theorem}

The above theorem is useful in verifying condition (\ref{cond1}). Indeed,
let the initial value $f\geq 0$ be chosen so that 
\begin{equation}
f(y)\leq KS_{\eta }\psi (y),\quad y\in D,  \label{K}
\end{equation}
where $\eta \geq 1$ is fixed and $K>0$ is a sufficiently small constant to
be specified later on. Therefore $S_{t}f\leq KS_{t+\eta }\psi $, and for any 
$t>0$, 
\begin{eqnarray*}
S_{t}f(x) &\leq &K\int_{D}p_{t+\eta }(x,y)\psi (y)\,dy \\
&=&K\int_{D}e^{\lambda _{1}(t+\eta )}\frac{p_{t+\eta }(x,y)}{\psi (x)\psi (y)%
}e^{-\lambda _{1}(t+\eta )}\psi (x)\psi^2 (y)\,dy \\
&\leq &K\left( \sup_{x\in D}\psi (x)\right) ^{2}\int_{D}e^{\lambda
_{1}(t+\eta )}\sup_{x,y\in D}\frac{p_{t+\eta }(x,y)}{\psi (x)\psi (y)}
e^{-\lambda _{1}(t+\eta )}\psi (y)\,dy \\
&\leq &K\left( \sup_{x\in D}\psi (x)\right) ^{2}\int_{D}\left( 1+c(1\wedge
(t+\eta ))^{-(d+2)/2}e^{-(\lambda _{2}-\lambda _{1})(t+\eta )}\right)
e^{-\lambda _{1}(t+\eta )}\psi (y)\,dy \\
&=&K\left( \sup_{x\in D}\psi (x)\right) ^{2}\left( e^{-\lambda _{1}(t+\eta
)}+ce^{-\lambda _{2}(t+\eta )}\right) \int_{D}\psi (y)\,dy \\
&\leq &K(1+c)e^{-\lambda _{1}\eta }\left( \sup_{x\in D}\psi (x)\right)
^{2}e^{-\lambda _{1}t}\int_{D}\psi (y)\,dy,
\end{eqnarray*}
which is independent of $x.$ Since the function $(t,x)\mapsto S_{t}f(x)$ is
uniformly bounded in $x$, condition (\ref{cond1}) is satisfied provided that 
\begin{equation*}
\Lambda \beta \left[ K(1+c)e^{-\lambda _{1}\eta }\left( \sup_{x\in D}\psi
(x)\right) ^{2}\int_{D}\psi (y)\,dy\right]^{\beta }\int_{0}^{\infty
}dr\,e^{\kappa \beta W_{r}-(\lambda _{1}+\kappa ^{2}/2)\beta r}<1,
\end{equation*}
or 
\begin{equation}
\int_{0}^{\infty }dr\,e^{\kappa \beta W_{r}-(\lambda _{1}+\kappa
^{2}/2)\beta r}\ < \ \frac{e^{\lambda _{1}\beta \eta }}{\displaystyle\Lambda
\beta \left[ K(1+c)\left( \sup_{x\in D}\psi (x)\right) ^{2}\int_{D}\psi
(y)\,dy\right]^{\beta }}.  \label{cond3}
\end{equation}
Notice that condition (\ref{cond2}) is satisfied if $K$ in (\ref{K}) is
sufficiently small.

We have proved the following

\begin{theorem}
\label{theorem5} Let $G$ satisfy (\ref{G2}), and let $D$ be a connected,
bounded $C^{1,\alpha }$-domain in $%
\mathbb{R}
^{d}$, where $\alpha >0.$ If (\ref{K}) and (\ref{cond3}) hold for some $%
\eta>0$ and $K>0$, then the solution of Equation (\ref{mild}) is global.
\end{theorem}

\noindent\textbf{Remarks.} 1) The integral on the left of (\ref{cond3})
coincides with the corresponding integral in sections 2 and 3. The same type
of bounds as in Section 3 can therefore be applied to estimate the
probability of existence of a global positive solution. By means of the law
of the iterated logarithm for $W$ we see from (\ref{cond3}) that the
presence of a noise may help to prevent blowup in finite time.

2) If $G(z)=\Lambda z^{1+\beta }$, the results of this section can be
applied to the solution $u$ of equation (\ref{E1}) since $v(t,x)=e^{-\kappa
W_{t}}u(t,x)$, $t\geq 0$, $x\in D.$ \bigskip

{\noindent\textbf{Acknowledgement}} \ M. Dozzi would like to thank Prof. Z.
Brze\'{z}niak for discussions on blowup of PDE's and SPDE's and acknowledges
the European Commission for partial support by the grant PIRSES 230804
"Multifractionality". J. A. L\'{o}pez-Mimbela acknowledges the hospitality
of \ Institut Elie Cartan, Nancy University, where part of this work was
done.

\vskip.5in

\begin{minipage}{3in}
\lineskip .25cm
\lineskiplimit .25cm

Marco Dozzi \\
IECN, Nancy Universit\'{e}s,\\
B.P. 239,\\
54506 Vandoeuvre-l\`{e}s-Nancy, France\\
{\tt dozzi@iecn.u-nancy.fr}
\end{minipage}
\begin{minipage}{3in}
\lineskip .25cm
\lineskiplimit .25cm

Jos\'{e} Alfredo L\'{o}pez-Mimbela\\
Centro de Investigaci\'{o}n en Matem\'{a}ticas\\
Apartado Postal 402\\ 36000 Guanajuato, Mexico\\
{\tt jalfredo@cimat.mx}
\end{minipage}


\begin{thebibliography}{99}
\bibitem{Bandle} Bandle, C.; Dozzi, M.; Schott, R. Blow up behaviour of a
stochastic partial differential equation of reaction-diffusion type.
Stochastic analysis: random fields and measure-valued processes (Ramat Gan,
1993/1995), 27--33, \emph{Israel Math. Conf. Proc.}, \textbf{10}, Bar-Ilan
Univ., Ramat Gan, 1996.

\bibitem{Berge} Berg\'{e}, B., Chueshov, I.D., Vuillermot, P.-A. On the
behaviour of solutions to certain parabolic SPDE's driven by Wiener
processes. \emph{Stochastic Process. Appl.} \textbf{92} (2001), 237-263.

\bibitem{Borodin and Salminen (2002)}  Borodin, Andrei N.; Salminen, Paavo 
\textsl{Handbook of Brownian motion---facts and  formulae.}  Second edition.
Probability and its Applications. Birkh\"{a}user Verlag, Basel, 2002.

\bibitem{Chueshov and Vuillermot} Chueshov, Igor D.; Vuillermot, Pierre A.
Long-time behavior of solutions to a class of stochastic parabolic equations
with  homogeneous white noise: It\^{o}'s case. \emph{Stochastic Anal. Appl.} 
\textbf{18} (2000), no. 4,  581--615.

\bibitem{Denis et al} Denis, Laurent; Matoussi, Anis; Stoica, Lucretiu. $L%
\sp p$ estimates for the uniform norm of solutions of quasilinear SPDE's. 
\emph{Probab. Theory Related Fields} \textbf{133} (2005), no. 4, 437--463.

\bibitem{Dufresne} Dufresne, Daniel. The distribution of a perpetuity, with
applications to risk theory and pension funding. \emph{Scand. Actuar. J.} 
\textbf{1990}, no. 1-2, 39--79.

\bibitem{Friedman} Friedman, Avner. \textsl{Partial differential equations
of parabolic type}. Prentice-Hall, Inc., Englewood Cliffs, N.J. 1964.

\bibitem{Fujita} Fujita, H. On the blowing up of solutions of the Cauchy
problem for $u_t=\Delta u + u^{1+\alpha}$. \emph{J. Fac. Sci. Univ. Tokyo,
Sect. 1} \textbf{13} (1966), 109--124.

\bibitem{Gyongy and Rovira} Gy\"{o}ngy, Istv\'{a}n; Rovira, Carles. On $L\sp %
p$-solutions of semilinear stochastic partial differential equations. \emph{%
Stochastic Process. Appl.} \textbf{90} (2000), no. 1, 83--108.

\bibitem{Klebaner} Klebaner, Fima C. \textsl{Introduction to stochastic
calculus with applications}. Second edition. Imperial College Press, London,
2005.

\bibitem{Krylov} Krylov, N. V. An analytic approach to SPDEs. Stochastic
partial differential equations: six perspectives, 185--242, \emph{Math.
Surveys Monogr.}, \textbf{64}, Amer. Math. Soc., Providence, RI, 1999.

\bibitem{Lototski
and Rozovskii} Lototsky, Sergey; Rozovskii, Boris.
Stochastic differential equations: a Wiener chaos approach. \textsl{From
stochastic calculus to mathematical finance}, 433--506, Springer, Berlin,
2006.

\bibitem{Manthey} Manthey, R., Zausinger, T. Stochastic evolution equations
in $L_{\rho}^{2\nu}$. \emph{Stoch. Stoch. Reports} \textbf{66} (1999),
37--85.

\bibitem{Mikulevicius and Pragarauskas} Mikulevicius, R.; Pragarauskas, H.
On Cauchy-Dirichlet problem for parabolic quasilinear SPDEs. \emph{Potential
Anal.} \textbf{25} (2006), no. 1, 37--75.

\bibitem{Ouhabaz and Wang} Ouhabaz, El Maati, Wang, Feng-Yu. Sharp estimates
for intrinsic ultracontractivity on $C\sp {1,\alpha}$-domains. \emph{%
Manuscripta Math.} \textbf{122} (2007), no. 2, 229--244.

\bibitem{Pazy} Pazy, A. \textsl{Semigroups of linear operators and
applications to partial differential equations}. Applied Mathematical
Sciences, 44. Springer-Verlag, New York, 1983.

\bibitem{Rothe} Rothe, F. \textsl{Global solutions of reaction-diffusion
systems}. Lecture Notes in Mathematics, 1072. Springer-Verlag, Berlin, 1984.

\bibitem{Yor} Yor, Marc. Exponential functionals of Brownian motion and
related processes.  Springer Finance. Springer-Verlag, Berlin, 2001.
\end{thebibliography}
\end{document}